\documentclass{article}
\begin{document}
\title{\bf 
On the intelligibility of the universe and the notions of 
simplicity, complexity and irreducibility
}

\author{Gregory Chaitin, \it IBM Research Division \\
\it chaitin@us.ibm.com}
\date{}
\maketitle
\begin{abstract}
We discuss views about whether the universe can be rationally comprehended,
starting with Plato, then Leibniz, and then the views of
some distinguished scientists of the previous century.
Based on this,
we defend the thesis that comprehension is compression, i.e., 
explaining many facts using few theoretical assumptions,
and that a theory may be viewed as a computer program for calculating observations.
This provides motivation for defining the complexity of something to be
the size of the simplest theory for it, in other words, the size of 
the smallest program for calculating it. This is the central
idea of algorithmic information theory (AIT), a field of theoretical computer science.
Using the mathematical concept of program-size complexity, we exhibit irreducible
mathematical facts, mathematical facts that cannot be demonstrated
using any mathematical theory simpler than they are.
It follows that
the world of mathematical ideas has infinite complexity and is therefore not
fully comprehensible, at least not in a static fashion.
Whether the physical world has finite or infinite complexity remains to
be seen. Current science believes that the world contains randomness, and is
therefore also infinitely complex, but a deterministic universe
that simulates randomness via pseudo-randomness is also a possibility,
at least according to recent highly speculative work of S. Wolfram. 
[Written for a meeting of the German Philosophical Society, Bonn, September 2002.]
\end{abstract}

\footnotesize
``Nature uses only the longest threads to weave her patterns,
so that each small piece of her fabric reveals the organization of the
entire tapestry.''
   
---Feynman, {\it The Character of Physical Law,} 1965,
at the very end of Chapter 1, ``The Law of Gravitation''.\footnote
{An updated version of this chapter would no doubt include a discussion
of the infamous astronomical missing mass problem.}
   
``The most incomprehensible thing about the universe
is that it is comprehensible.''
    
---Attributed to Einstein. 
The original source, where the wording is somewhat different, is Einstein, ``Physics and Reality'',
1936, reprinted in Einstein, {\it Ideas and Opinions,} 1954.\footnote
{Einstein actually wrote
   ``Das ewig Unbegreifliche an der Welt ist ihre Begreiflichkeit''.
Translated word for word, this is
   ``The eternally incomprehensible about the world is its comprehensibility''.
But I prefer the version given above, which emphasizes the paradox.}
\\
    
\normalsize
It's a great pleasure for me to speak at this meeting of the German Philosophical Society.
Perhaps it's not generally known that at the end of his life 
my predecessor Kurt G\"odel was obsessed with Leibniz.\footnote
{See Menger, {\it Reminiscences of the Vienna Circle and the Mathematical Colloquium,} 1994.}
Writing this paper was for me a voyage of discovery---of the depth of
Leibniz's thought!
Leibniz's power as a philosopher is informed by his genius as a mathematician; 
as I'll explain, some of the key ideas of AIT are clearly
visible in embryonic form in his 1686 {\it Discourse on Metaphysics.}

\section*{\bf
I Plato's {\it Timaeus}---The Universe is Intelligible. Origins of the Notion of
Simplicity: Simplicity as Symmetry [Brisson, Meyerstein 1991]
}

\footnotesize
``[T]his is the central idea developed in the {\it Timaeus\/}: the order
established by the demiurge in the universe becomes manifest as the symmetry
found at its most fundamental level, a symmetry which makes possible a
mathematical description of such a universe.''
    
---Brisson, Meyerstein, {\it Inventing the Universe,} 1995 (1991 in French).
This book discusses
the cosmology of
Plato's {\it Timaeus,} modern cosmology and AIT; one of their key insights
is to identify symmetry with simplicity.
\\

\normalsize
According to 
Plato, the world is rationally understandable
because 
it has structure.
And the universe has structure, because it is a work of art
created by a God who is a mathematician. 
Or, more abstractly, the structure of the world consists of God's thoughts,
which are mathematical. The fabric of reality is built out of eternal mathematical truth.
[Brisson, Meyerstein, {\it Inventer
l'Univers,} 1991]
    
{\it Timaeus\/} postulates that simple, symmetrical geometrical forms
are the building blocks for the universe:
the circle and 
the regular solids (cube, tetrahedron, icosahedron, dodecahedron).
    
What was the evidence 
that convinced the ancient Greeks that the world is comprehensible?
Partly it was the beauty of mathematics,
particularly geometry and number theory,
and partly the Pythagorean work on the physics of stringed instruments
and musical tones, and in astronomy, the regularities in the motions of the planets and
the starry heavens and eclipses.
Strangely enough, mineral crystals, whose symmetries magnify enormously 
quantum-mechanical
symmetries that are found at the
atomic and molecular level, are never mentioned.
   
What is our current cosmology?
   
Since the chaos of everyday existence provides little evidence of simplicity,
biology is based
on chemistry is based on physics is based on high-energy or particle
physics. The attempt to
find underlying simplicity and pattern leads reductionist modern science
to break things into
smaller and smaller components in an effort to find the underlying simple
building blocks.
   
And the modern version of the cosmology of {\it Timaeus\/} 
is the application of symmetries or group theory
to understand sub-atomic particles
(formerly called elementary particles), for example, Gell-Mann's
eightfold way, which
predicted new particles. 
This work classifying the ``particle zoo'' also resembles
Mendeleev's
periodic table of the elements that organizes their
chemical properties so well.\footnote
{For more on this,
see the essay by Freeman Dyson on ``Mathematics in the Physical Sciences''
in COSRIMS, {\it The Mathematical Sciences,} 1969.
This is an article of his that was originally published in {\it Scientific American.}}
   
And modern physicists have also come up with a possible answer to
the Einstein quotation at the beginning of this paper.  
Why do they think that the universe is comprehensible?
They invoke the so-called ``anthropic principle''
[Barrow, Tipler, {\it The Anthropic Cosmological Principle,} 1986], 
and declare that we would not be here to
ask this question unless the universe had enough order for complicated creatures
like us to evolve!
   
Now let's proceed to the next major step in the evolution of ideas on simplicity and
complexity,
which is a stronger version of the Platonic creed due to Leibniz.

\section*{\bf
II What Does it Mean for the Universe to be Intelligible? Leibniz's 
Discussion of Simplicity, Complexity and Lawlessness [Weyl 1932]
}

\footnotesize
``As for the simplicity of the ways of God, this holds
properly with respect to his means, as opposed to the variety, richness,
and abundance, which holds with respect to his ends or effects.''
    
``But, {\bf when a rule is extremely complex, what is in
conformity with it passes for irregular}. Thus, one can say, in whatever
manner God might have created the world, it would always have been regular
and in accordance with a certain general order. But {\bf God has chosen the
most perfect world, that is, the one which is at the same time the simplest
in hypotheses and the richest in phenomena}, as might be a line in geometry
whose construction is easy and whose properties and effects are extremely
remarkable and widespread.''
   
---Leibniz, {\it Discourse on Metaphysics,} 1686, 
Sections 5--6, from
Leibniz, {\it Philosophical Essays,} edited and translated by Ariew and
Garber, 1989, pp.\ 38--39.
   
``The assertion that nature is governed by strict laws
is devoid of all content if we do not add the statement that it is governed
by mathematically simple laws\ldots\ That {\bf the notion of law becomes empty
when an arbitrary complication is permitted} was already pointed out
by Leibniz in his {\it Metaphysical Treatise\/} [{\it Discourse on Metaphysics\/}].
Thus simplicity becomes a working principle in the natural sciences\ldots\
The astonishing thing is not that there exist natural laws, but that the
further the analysis proceeds, the finer the details, the finer the elements
to which the phenomena are reduced, the simpler---and not the more complicated,
as one would originally expect---the fundamental relations become and the
more exactly do they describe the actual occurrences. But this circumstance
is apt to weaken the metaphysical power of determinism, since it makes
the meaning of natural law depend on the fluctuating distinction between
mathematically simple and complicated functions or classes of functions.''
    
---Hermann Weyl, {\it The Open World, Three Lectures on
the Metaphysical Implications of Science,} 1932, pp.\ 40--42. See a similar
discussion on pp.\ 190--191 of Weyl, {\it Philosophy of Mathematics and Natural
Science,} 1949, Section 23A, ``Causality and Law''.\footnote{This is a remarkable
anticipation of my definition of ``algorithmic randomness'', as a set of observations
that only has what Weyl considers to be unacceptable theories, ones that are 
as complicated as the observations themselves, without any ``compression''.}
    
``Weyl said, not long ago, that `the problem of simplicity
is of central importance for the epistemology of the natural sciences'.
Yet it seems that interest in the problem has lately declined; perhaps
because, especially after Weyl's penetrating analysis, there seemed to
be so little chance of solving it.''
    
---Weyl, {\it Philosophy of Mathematics and Natural Science,}
1949, p.\ 155, quoted in Popper, {\it The Logic of Scientific Discovery,}
1959, Chapter VII, ``Simplicity'', p.\ 136.
\\

\normalsize
In his novel {\it Candide,} Voltaire ridiculed Leibniz, caricaturing Leibniz's subtle
views with the memorable phrase
{\bf ``this is the best of all possible worlds''}.
Voltaire also ridiculed the efforts of Maupertius
to develop a physics in line with Leibniz's views, one based on a principle
of least effort.
   
Nevertheless versions of least effort play a fundamental
role in modern science, starting with Fermat's deduction of the laws for
reflection and refraction of light from a principle of least time.
This
continues with the Lagrangian formulation of mechanics, stating that the actual
motion minimizes the integral of the difference between the potential and
the kinetic energy. And least effort is even important at the current frontiers,
such 
as in Feynman's path integral formulation of quantum mechanics (electron
waves) and quantum electrodynamics (photons, electromagnetic field
quanta).\footnote 
{See the short discussion of minimum principles in Feynman, {\it The Character
of Physical Law,} 1965, Chapter 2, ``The Relation of Mathematics to Physics''. For
more information, see
{\it The Feynman Lectures on Physics,} 1963,
Vol.\ 1, Chapter 26, ``Optics: The Principle of Least Time'',
Vol.\ 2, Chapter 19, ``The Principle of Least Action''.}
   
However, all this modern physics refers to versions of least effort, not to ideas, not to
information, and not to complexity---which are more closely connected with Plato's
original emphasis on symmetry and intellectual simplicity = intelligibility.
An analogous situation occurs in theoretical computer science,
where work on computational complexity is usually focussed
on time, not on the complexity of ideas or information. Work on time complexity
is of great practical value, but I believe that the complexity of ideas is of
greater conceptual significance. 
Yet another example of the effort/information divide is the fact that 
I am interested in the irreducibility
of ideas (see Sections V and VI),  
while Stephen Wolfram (who is discussed later in this section) 
instead emphasizes time irreducibility, physical systems
for which there are no predictive short-cuts
and the fastest way to see what they do is just to run them. 
    
Leibniz's doctrine concerns more than ``least effort'',
it also implies that the ideas that produce
or govern this world are as beautiful and as simple as possible.
In more modern terms, God employed the smallest possible amount of intellectual material to build the world,
and the laws of physics are as simple and as beautiful as they can be and
allow us, intelligent beings, to evolve.\footnote{This is a kind of ``anthropic principle'',
the attempt to deduce things about the universe from the fact that we are
here and able to look at it.} The belief in this Leibnizean doctrine lies behind 
the continuing reductionist
efforts of high-energy physics (particle physics) to find the ultimate components of
reality. 
The continuing vitality of this Leibnizean doctrine
also lies behind astrophysicist
John Barrow's emphasis in his ``Theories of Everything'' essay on finding the minimal TOE that explains
the universe, a TOE that is as simple as possible, with no redundant elements
(see Section VII below).
    
{\bf Important point:} To say that the fundamental laws of physics
must be simple does not at all imply that it is easy or fast to deduce from
them how the world works, that it is quick to make predictions from the
basic laws. The {\it apparent complexity\/} of the world we live in---a 
phrase that is constantly repeated in Wolfram, {\it A New Kind of Science,} 
2002---then comes from the long deductive path from the basic laws to the level
of our experience.\footnote{It could also come from the complexity of the initial conditions,
or from coin-tossing, i.e., randomness.} So again, I claim that minimum information is
more important than minimum time, which is why in Section IV I do not care how long a
minimum-size program takes to produce its output, 
nor how much time it takes to calculate experimental data using a scientific theory.
   
{\bf More on Wolfram:}  In 
{\it A New Kind of Science,} Wolfram reports on his systematic computer
search for simple rules with very complicated consequences, very much in the
spirit of Leibniz's remarks above.  
First Wolfram amends the Pythagorean insight that Number rules the universe to assert the primacy
of Algorithm, not Number. 
And those are {\bf discrete} algorithms, it's a digital philosophy!\footnote
{That's a term invented by Edward Fredkin, who has worked on related ideas.}
Then Wolfram sets out to survey all possible worlds, at least all the simple ones.\footnote
{That's why his book is so thick!}
Along the way he finds a lot of interesting stuff.
For example,
Wolfram's cellular automata rule 110 is a universal computer, 
an amazingly simple one,
that can carry out {\bf any} computation.
{\it A New Kind of Science\/} is an attempt to discover the laws of the 
universe by pure thought, to search systematically for God's building blocks!
    
{\bf The limits of reductionism:} In what sense can biology and psychology
be reduced to mathematics and physics?! This is indeed the acid test of a
reductionist viewpoint! Historical contingency is often invoked here: life
as ``frozen accidents'' (mutations), not something fundamental [Wolfram, Gould]. 
Work on artificial life (Alife) plus
advances in robotics are particularly aggressive reductionist attempts. The
normal way to ``explain'' life is evolution by natural selection, ignoring
Darwin's own sexual selection and symbiotic/cooperative views of the origin
of biological progress---new species---notably espoused by Lynn Margulis (``symbiogenesis''). 
Other problems
with Darwinian gradualism: following the DNA as software paradigm,
small changes in DNA software can produce big changes in organisms, and a good way to build
this software is by trading useful subroutines (this is called horizontal or lateral
DNA transfer).\footnote{This is how bacteria acquire immunity to antibiotics.}
In fact, there is a lack of fossil evidence
for many intermediate forms,\footnote{Already noted by Darwin.} 
which is evidence for rapid production of new species 
(so-called ``punctuated equilibrium'').

\section*{\bf
III What do Working Scientists Think about Simplicity and Complexity?
}

\footnotesize
``Science itself, therefore, may be regarded as a minimal
problem, consisting of the completest possible presentment of facts with
the {\it least possible expenditure of thought\ldots\ } Those ideas that hold
good throughout the widest domains of research and that supplement the
greatest amount of experience, are {\it the most scientific.}''
    
---Ernst Mach, {\it The Science of Mechanics,} 1893,
Chapter IV, Section IV, ``The Economy of Science'', 
reprinted in Newman, {\it The World of Mathematics,}
1956.
    
``Furthermore, the attitude that theoretical physics does
not explain phenomena, but only classifies and correlates, is today accepted
by most theoretical physicists. This means that the criterion of success
for such a theory is simply whether it can, by a simple and elegant classifying
and correlating scheme, cover very many phenomena, which without this scheme
would seem complicated and heterogeneous, and whether the scheme even covers
phenomena which were not considered or even not known at the time when
the scheme was evolved. (These two latter statements express, of course,
the unifying and the predicting power of a theory.)''
   
---John von Neumann, ``The Mathematician'', 1947,
reprinted in Newman, {\it The World of Mathematics,} 1956, and in 
Br\'ody, V\'amos, {\it The Neumann Compendium,} 1995.
   
``{\bf These fundamental concepts and postulates, which cannot
be further reduced logically, form the essential part of a theory, which
reason cannot touch. It is the grand object of all theory to make these
irreducible elements as simple and as few in number as possible}\ldots\ [As]
the distance in thought between the fundamental concepts and laws on the
one side and, on the other, the conclusions which have to be brought into
relation with our experience grows larger and larger, the simpler the logical
structure becomes---that is to say, the smaller the number of logically
independent conceptual elements which are found necessary to support the
structure.''
   
---Einstein, ``On the Method of Theoretical Physics'', 1934,
reprinted in Einstein, {\it Ideas and Opinions,} 1954.
    
``The aim of science is, on the one hand, a comprehension,
as {\it complete\/} as possible, of the connection between the sense experiences
in their totality, and, on the other hand, the accomplishment of this aim
{\it by the use of a minimum of primary concepts and relations.} (Seeking as
far as possible, logical unity in the world picture, i.e., paucity in logical
elements.)''
    
``Physics constitutes a logical system of thought which
is in a state of evolution, whose basis cannot be distilled, as it were,
from experience by an inductive method, but can only be arrived at by free
invention\ldots\ Evolution is proceeding in the direction of increased simplicity
of the logical basis. In order further to approach this goal, we must resign
to the fact that the logical basis departs more and more from the facts
of experience, and that the path of our thought from the fundamental basis
to those derived propositions, which correlate with sense experiences,
becomes continually harder and longer.''
    
---Einstein, ``Physics and Reality'', 1936, reprinted in
Einstein, {\it Ideas and Opinions,} 1954.
   
``[S]omething general will have to be said\ldots\ about the
points of view from which physical theories may be analyzed critically\ldots\
The first point of view is obvious: the theory must not contradict empirical
facts\ldots\ The second point of view is not concerned with the relationship
to the observations but with the premises of the theory itself, with what
may briefly but vaguely be characterized as the `naturalness' or `logical
simplicity' of the premises (the basic concepts and the relations between
these)\ldots\ We prize a theory more highly if, from the logical standpoint,
it does not involve an arbitrary choice among theories that are equivalent
and possess analogous structures\ldots\ I must confess herewith that I cannot
at this point, and perhaps not at all, replace these hints by more precise
definitions. I believe, however, that a sharper formulation would be possible.''
   
---Einstein, ``Autobiographical Notes'', originally published in Schilpp,
{\it Albert Einstein, Philosopher-Scientist,} 1949, and reprinted
as a separate book in 1979.
    
``What, then, impels us to devise theory after theory? Why do we devise
theories at all? The answer to the latter question is simply: because {\bf we enjoy
`comprehending,' i.e., reducing phenomena by the process of logic to something
already known} or (apparently) evident. New theories are first of all necessary
when we encounter new facts which cannot be `explained' by existing theories.
But this motivation for setting up new theories is, so to speak, trivial, imposed
from without. There is another, more subtle motive of no less importance. This is
the striving toward unification and simplification of the premises of the theory
as a whole (i.e., Mach's principle of economy, interpreted as a logical principle).''
    
``There exists a passion for comprehension, just as there exists a passion for music.
That passion is rather common in children, but gets lost in most people later on.
Without this passion, there would be neither mathematics nor natural science.
Time and again the passion for understanding has led to the illusion that man is
able to comprehend the objective world rationally, by pure thought, without any
empirical foundations---in short, by metaphysics. I believe that every true theorist
is a kind of tamed metaphysicist, no matter how pure a `positivist' he may fancy 
himself. {\bf The metaphysicist believes that the logically simple is also the real.
The tamed metaphysicist believes that not all that is logically simple is embodied
in experienced reality, but that the totality of all sensory experience can be
`comprehended' on the basis of a conceptual system built on premises of great
simplicity.} The skeptic will say that this is a `miracle creed.' Admittedly so, but 
it is a miracle creed which has been borne out to an amazing extent by the 
development of science."
   
---Einstein, ``On the Generalized Theory of Gravitation'', 1950, reprinted in 
Einstein, {\it Ideas and Opinions,} 1954.
    
``One of the most important things in this `guess---compute consequences---compare
with experiment' business is to know when you are right. It is possible to know
when you are right way ahead of checking all the consequences. You can recognize
truth by its beauty and simplicity. It is always easy when you have made a guess,
and done two or three little calculations to make sure that it is not obviously
wrong, to know that it is right. {\bf When you get it right, it is obvious that it is
right---at least if you have any experience---because usually what happens is that
more comes out than goes in.} Your guess is, in fact, that something is very
simple. If you cannot see immediately that it is wrong, and it is simpler than it 
was before, then it is right. The inexperienced, and crackpots, and people like
that, make guesses that are simple, but you can immediately see that they are wrong,
so that does not count. Others, the inexperienced students, make guesses that are
very complicated, and it sort of looks as if it is all right, but I know it is not
true because the truth always turns out to be simpler than you thought. What we
need is imagination, but imagination in a terrible strait-jacket. We have to find
a new view of the world that has to agree with everything that is known, but
disagree in its predictions somewhere, otherwise it is not interesting. And in that
disagreement it must agree with nature\ldots''
    
---Feynman, {\it The Character of Physical Law,} 1965, Chapter 7, ``Seeking
New Laws''.
    
``It is natural that a man should consider the work of his hands or his brain to
be useful and important. Therefore nobody will object to an ardent experimentalist
boasting of his measurements and rather looking down on the `paper and ink' physics
of his theoretical friend, who on his part is proud of his lofty ideas and despises
the dirty fingers of the other. But in recent years this kind of friendly rivalry
has changed into something more serious\ldots\ [A] school of extreme experimentalists\ldots\
has gone so far as to reject theory altogether\ldots\ There is also a movement in the
opposite direction\ldots\ claiming that to the mind well trained in mathematics and
epistemology the laws of Nature are manifest without appeal to experiment.''
    
``Given the knowledge and the penetrating brain of our mathematician, Maxwell's
equations are a result of pure thinking and the toil of experimenters antiquated
and superfluous. I need hardly explain to you the fallacy of this standpoint.
It lies in the fact that none of the notions used by the mathematicians, such as
potential, vector potential, field vectors, Lorentz transformations, quite apart
from the principle of action itself, are evident or given {\it a priori.}
Even if an extremely gifted mathematician had constructed them to describe the
properties of a possible world, neither he nor anybody else would have had the
slightest idea how to apply them to the real world.''
   
``Charles Darwin, my predecessor in my Edinburgh chair, once said something like
this: `The Ordinary Man can see a thing an inch in front of his nose; a few
can see things 2 inches distant; if anyone can see it at 3 inches, he is a
man of genius.' I have tried to describe to you some of the acts of these 2-
or 3-inch men. My admiration of them is not diminished by the consciousness
of the fact that they were guided by the experience of the whole human race to
the right place into which to poke their noses. {\bf I have also not endeavoured
to analyse the idea of beauty or perfection or simplicity of a natural law}
which has often guided the correct divination. {\bf I am convinced that such an
analysis would lead to nothing; for these ideas are themselves subject to
development. We learn something new from every new case}, and I am not inclined
to accept final theories about invariable laws of the human mind.''
    
``My advice to those who wish to learn the art of scientific prophecy is not
to rely on abstract reason, but to decipher the secret language of Nature from
Nature's documents, the facts of experience.''
   
---Max Born, {\it Experiment and Theory in Physics,} 1943, pp.\ 1, 8, 34--35, 44. 
\\

\normalsize
These eloquent discussions of the role that simplicity and complexity
play in scientific discovery by these
distinguished 20th century scientists show the importance that they ascribe to these 
questions.
    
In my opinion, the
fundamental point is this: The belief that the universe is rational,
lawful, is of no value if the laws are
too complicated for us to comprehend, and is even meaningless if the
laws are as complicated
as our observations, since the laws are then no simpler than the world
they are supposed to explain.
As we saw in the previous section,
this was emphasized (and attributed to Leibniz) by Hermann Weyl, a fine
mathematician and mathematical
physicist.
    
But perhaps we are overemphasizing the role that the notions of simplicity
and complexity play in science?
     
In his beautiful 1943 lecture published as a small book
on {\it Experiment and Theory in Physics,}
the theoretical physicist Max Born criticized those who think that we can
understand Nature by pure thought, without hints from experiments.  In particular,
he was referring to now forgotten and rather fanciful theories put forth by
Eddington and Milne.  Now he might level these criticisms at string theory and 
at Stephen Wolfram's {\it A New Kind of Science\/}
[Jacob T. Schwartz, private communication].
    
Born has a point.
Perhaps the universe {\bf is} complicated, not simple! This certainly
seems to be the case in biology more than in physics.  
Then thought alone is insufficient; we need empirical data.
But simplicity certainly
reflects what we mean by understanding: {\bf understanding is compression}.
So perhaps this is more about the human mind than it is about the universe.
Perhaps our emphasis on simplicity says more about us than it says about the universe!
    
Now we'll try to capture some of the essential features of these philosophical
ideas in a mathematical theory.

\section*{\bf
IV A Mathematical Theory of Simplicity, Complexity and Irreducibility: AIT
}

The basic idea of algorithmic information theory (AIT) is that a scientific theory
is a computer program, and the smaller, the more concise the program is, the better the theory!
    
But the idea is actually much broader than that.
{\bf The central idea of algorithmic information theory is reflected in the belief
that the following diagrams all have something fundamental in common.}
In each case, ask how much information we put in versus how much we get out.
And everything is digital, discrete.
    
Shannon information theory (communications engineering), noiseless coding:
\begin{center}
   encoded message $\rightarrow$ {\bf Decoder} $\rightarrow$ original message
\end{center}
Model of scientific method:
\begin{center}
   scientific theory $\rightarrow$ {\bf Calculations} $\rightarrow$ empirical/experimental data
\end{center}
Algorithmic information theory (AIT), definition of program-size complexity:
\begin{center}
   program $\rightarrow$ {\bf Computer} $\rightarrow$ output
\end{center}
Central dogma of molecular biology: 
\begin{center}
   DNA $\rightarrow$ {\bf Embryogenesis/Development} $\rightarrow$ organism
\end{center}
(In this connection, see K\"uppers, {\it Information and the Origin of Life,} 1990.) 
Turing/Post abstract formulation of a Hilbert-style formal axiomatic mathematical theory
as a mechanical procedure for systematically deducing all possible consequences from the axioms:
\begin{center}
   axioms $\rightarrow$ {\bf Deduction} $\rightarrow$ theorems
\end{center}
Contemporary physicists' efforts to find a Theory of Everything (TOE):
\begin{center}
   TOE $\rightarrow$ {\bf Calculations} $\rightarrow$ Universe
\end{center}
Leibniz, {\it Discourse on Metaphysics,} 1686:
\begin{center}
   Ideas $\rightarrow$ {\bf Mind of God} $\rightarrow$ The World
\end{center}
In each case the left-hand side is smaller, much smaller, than the right-hand
side. In each case, the right-hand side can be constructed (re-constructed)
mechanically, or systematically, from the left-hand side. And in each case
we want to keep the right-hand side fixed while making the left-hand side
as small as possible.
Once this is accomplished,
we can use the size of the left-hand side as 
a measure of the simplicity or the complexity
of the corresponding right-hand side.
    
Starting with this one simple idea, 
of looking at the size of computer programs, or at program-size complexity,
you can develop a sophisticated, elegant
mathematical theory, AIT, as you can see in my four Springer-Verlag volumes
listed in the bibliography of this paper.
   
But, I must confess that AIT makes a large number of {\bf important hidden assumptions!}
What are they?
    
Well, one
important hidden assumption of AIT is that the choice of computer or of
computer programming language is not too important, that it does not
affect program-size complexity too much, in any fundamental way.
This is debatable.
    
Another important tacit assumption: we use the discrete computation approach of
Turing 1936, eschewing computations with ``real'' (infinite-precision) numbers
like $\pi$ = 3.1415926\ldots\ which have an infinite number of digits when
written in decimal notation, but which correspond, from a geometrical point of
view, to a single point on a line, an elemental notion in continuous, but not in
discrete, mathematics.  Is the universe {\bf discrete} or {\bf continuous}?  
Leibniz is famous for his work on continuous mathematics.
AIT sides with the discrete, not with the continuous.
[Fran\c{c}oise Chaitin-Chatelin, private communication]
    
Also, in AIT we completely ignore the {\bf time} taken by a computation,
concentrating only on the {\bf size} of the program.  And the computation run-times may be
monstrously large, quite impracticably so, in fact, totally astronomical in size.  
But trying to take time into account
destroys AIT, an elegant, simple theory of complexity, 
and one which imparts much intuitive understanding. So I think that it is
a mistake to try to take time into account when thinking about this kind of complexity.
    
We've talked about simplicity and complexity, but what about {\bf irreducibility}?
Now let's apply AIT to mathematical logic and obtain some limitative metatheorems.
However, following Turing 1936 and Post 1944, I'll use
the notion of algorithm
to deduce limits to formal reasoning, not G\"odel's original 1931 approach.
I'll take the position that a Hilbert-style mathematical theory, a
formal axiomatic theory, is a mechanical procedure for systematically
generating all the theorems by running through all possible proofs,
systematically deducing all consequences of the axioms.\footnote
{In a way, this point of view was anticipated by Leibniz with his
{\it lingua characteristica universalis.}}
Consider the size in bits of the algorithm for doing this. This is
how we measure the simplicity or complexity of the formal axiomatic theory. 
It's just another instance of program-size complexity!
   
But at this point, Chaitin-Chatelin insists, I should admit
that we are making an extremely embarrassing hidden assumption, which is
that you can systematically run through all the proofs.  
This assumption, which is bundled into my definition of a formal axiomatic theory,
means that we
are assuming that the language of our theory is static, and that no new
concepts can ever emerge.  But no human language or field of thought is static!\footnote
{And computer programming languages aren't static either, which can be quite a nuisance.}
And this idea of being able to make a numbered list with all possible proofs
was clearly anticipated by \'Emile Borel in 1927
when he pointed out that there is 
a real number 
with the problematical property that its
$N$th digit after the decimal point gives us the answer
to the $N$th yes/no question in French.\footnote
{Borel's work was brought to my attention
by Vladimir Tasi\'c in
his book {\it Mathematics and the Roots of Postmodern Thought,} 2001,
where he points out that in some ways it anticipates the $\Omega$
number that I'll discuss in Section IX.
Borel's paper is reprinted in
Mancosu, {\it From Brouwer to Hilbert,} 1998, pp.\ 296--300.}
    
Yes, I agree, a Hilbert-style formal axiomatic theory is indeed a fantasy,
but it is a fantasy that inspired many people, and one that even helped to lead
to the creation of modern programming languages.  It is a fantasy that it is
useful to take seriously long enough for us to show in Section VI that even
if you are willing to accept all these tacit assumptions, something else is terribly wrong.
Formal axiomatic theories can be criticized from within, as well as from without.
And it is far from clear how weakening these tacit assumptions would make it easier
to prove the irreducible mathematical truths that are exhibited in Section VI.
    
And the idea of a fixed, static computer programming language in which you write the
computer programs whose size you measure is also a fantasy.  Real computer programming
languages don't stand still, they evolve, and the size of the computer program
you need to perform a given task can therefore change.
Mathematical models of the world like these are always approximations, ``lies that
help us to see the truth'' (Picasso).  Nevertheless, if done properly, they can impart insight and
understanding, they can help us to comprehend, they can reveal unexpected connections\ldots

\section*{\bf
V From Computational Irreducibility to Logical Irreducibility.
Examples of Computational Irreducibility: ``Elegant'' Programs
}

Our goal in this section and the next is to use AIT to establish the existence
of {\bf irreducible mathematical truths}.  What are they, and why are they important?
    
Following Euclid's {\it Elements,} a mathematical
truth is established
by reducing it to simpler truths until self-evident truths---``axioms''
or ``postulates''\footnote{Atoms of thought!}---are reached. 
Here we exhibit an extremely large class of 
mathematical truths that are not at all self-evident but which
are {\bf not} consequences of any principles simpler than they are.
    
Irreducible truths are
highly problematical for traditional philosophies of mathematics,
but as discussed in Section VIII,
they can be
accommodated in an emerging ``quasi-empirical'' school of the foundations
of mathematics, which says that physics and mathematics are not that different. 
    
Our path to logical irreducibility 
starts with computational irreducibility.
Let's start by calling a computer program ``elegant'' if no smaller
program in the same language produces exactly the same output.
There are lots of elegant programs, at least one for each output. 
And it doesn't matter how {\bf slow} an elegant program is,
all that matters is that it be as {\bf small} as possible.
    
{\it 
An elegant program viewed as an object in its own right is computationally
irreducible. 
}
Why? Because otherwise you can get a more concise program for
its output by computing it first and then running it.
Look at this diagram:
\begin{center}
   program${}_2$ $\rightarrow$ {\bf Computer} $\rightarrow$ 
   program${}_1$ $\rightarrow$ {\bf Computer} $\rightarrow$ output
\end{center}
If program${}_1$ is as concise as possible, then program${}_2$
cannot be much more
concise than program${}_1$. Why?
Well, consider a fixed-sized routine for running a program and then immediately running its output.
Then
\begin{center}
   program${}_2$ + fixed-size routine
   $\rightarrow$ {\bf Computer} $\rightarrow$ output
\end{center}
produces exactly the same output as program${}_1$
and would be a more concise program for producing that output than program${}_1$ is.
But this is impossible because it contradicts our hypothesis that program${}_1$
was already as small as possible.
{\it Q.E.D.}
    
Why should elegant programs interest philosophers?
Well, because of
Occam's razor, because the best theory to explain a fixed set of data is an elegant program!
    
But how can we get irreducible truths?  Well, just try {\bf proving} that a program is elegant!

\section*{\bf
VI Irreducible Mathematical Truths.
Examples of Logical Irreducibility: Proving a Program is Elegant 
}

{\bf Hauptsatz:}
{\it 
You cannot prove that a program is elegant if its size is substantially
larger than the size of the algorithm for generating all the theorems in your
theory.
}
    
{\bf Proof:} 
The basic idea is to
run the first provably elegant program you encounter when you systematically
generate all the theorems, and that is substantially larger than the
size of the algorithm for generating all the theorems. Contradiction, unless no
such theorem can be demonstrated, or unless the theorem is false.

Now I'll explain why this works. We are given a formal axiomatic mathematical theory:
\begin{center}
   theory = program $\rightarrow$ {\bf Computer} $\rightarrow$ set of all theorems
\end{center}
We may suppose that this theory is an elegant program, 
i.e., as concise as possible for producing the set of theorems that it does.
Then the size of this program is by definition the complexity of the theory, since it is
the size of the smallest program for systematically generating the set of all the theorems,
which are all the consequences of the axioms.
Now consider a fixed-size routine with the property that
\begin{center}
   theory + fixed-size routine $\rightarrow$ {\bf Computer} $\rightarrow$
\\  
   output of the first provably elegant program larger than 
\\  
   complexity of theory
\end{center}
More precisely,
\begin{center}
   theory + fixed-size routine $\rightarrow$ {\bf Computer} $\rightarrow$
\\   
   output of the first provably elegant program larger than 
\\  
   (complexity of theory + size of the fixed-size routine)
\end{center}
This proves our assertion that a mathematical theory cannot prove that a program is elegant
if that program is substantially larger than the complexity of the theory.
    
Here is the proof of this result in more detail. 
The fixed-size routine knows its own size and is given the theory,
a computer program for generating theorems,
whose size it measures and which it then runs, until the first theorem is encountered
asserting that a particular program $P$ is elegant that is larger than the total input to
the computer.  
The fixed-size routine then runs the program $P$, and finally
produces as output the same output as $P$ produces.
But this is impossible, because the output from $P$ cannot be obtained from a program
that is smaller than $P$ is, not if, as we assume by hypothesis, 
all the theorems of the theory are true and $P$ is actually elegant.
Therefore $P$ cannot exist.
In other words,
if there is a provably elegant program $P$ whose
size is greater than the complexity of the theory + the size of
this fixed-size routine, either $P$ is actually inelegant or we have a contradiction.
{\it Q.E.D.}
    
Because no mathematical theory of finite complexity can enable you
to determine all the elegant programs,
the following is immediate:
    
{\bf Corollary:} {\it The mathematical universe has infinite complexity.}\footnote
{On the other hand, our current mathematical theories are {\bf not} very complex.
On pages 773--774 of {\it A New Kind of Science,} Wolfram makes this point by
exhibiting essentially all of the axioms for traditional mathematics---in just two pages!
However, a program to generate all the theorems would be larger.}
    
This strengthens G\"odel's 1931 refutation of Hilbert's belief that a single,
fixed formal axiomatic theory could capture all of mathematical truth.
    
Given the significance of this conclusion, it is natural to demand more information.
You'll notice that
I never said {\bf which} computer programming language I was using!
    
Well, you can actually carry out this proof using either
high-level languages such as the version of
LISP that I use in {\it The Unknowable,} or using low-level binary machine languages,
such as the one that I use in {\it The Limits of Mathematics.}
In the case of a high-level computer programming language, one measures the
size of a program in characters (or 8-bit bytes) of text.  In the case of a binary machine
language, one measures the size of a program in 0/1 bits.
My proof works either way.
    
But I must confess that
not all programming languages permit my proof to work out this neatly.  
The ones that do are the kinds
of programming languages that you use in AIT, the ones for which program-size complexity
has elegant properties instead of messy ones, the ones that directly expose the fundamental 
nature of this complexity concept
(which is also called algorithmic information content), 
not the programming languages that bury the basic idea in a mass of messy technical details.
    
This paper
started with philosophy, and then we developed a mathematical theory.
Now let's go back to philosophy.  In the last three sections of this paper 
we'll discuss the philosophical implications of AIT.

\section*{\bf
VII Could We Ever Be Sure that We Had the Ultimate TOE? [Barrow 1995]
}

\footnotesize
``The search for a `Theory of Everything' is the quest
for an ultimate compression of the world. Interestingly, Chaitin's proof
of G\"odel's incompleteness theorem using the concepts of complexity
and compression reveals that G\"odel's theorem is equivalent to the
fact that one cannot prove a sequence to be incompressible. We can never
prove a compression to be the ultimate one; there might be a yet deeper
and simpler unification waiting to be found.''
    
---John Barrow, essay on ``Theories of Everything'' in Cornwell, {\it Nature's
Imagination,} 1995, reprinted in Barrow, {\it Between Inner Space and
Outer Space,} 1999.
\\

\normalsize
Here is the first philosophical application of AIT.
According to astrophysicist
John Barrow, my work implies that even if we had
the optimum,
perfect, minimal (elegant!\@) TOE, we could never be sure a simpler theory would
not have the same
explanatory power.
     
(``Explanatory power'' is a pregnant phrase, and one can make a case that it is a better
name to use than the dangerous word ``complexity'', which has many other possible meanings.
One could then speak of a theory with $N$ bits of algorithmic explanatory
power, rather than describe it as a theory having a program-size complexity 
of $N$ bits. [Fran\c{c}oise Chaitin-Chatelin, private communication])
   
Well, you can dismiss Barrow by saying that the idea of having the ultimate TOE
is pretty crazy---who expects to be able to read the mind of God?!
Actually, Wolfram believes
that a systematic computer search might well find the ultimate TOE.\footnote  
{See pages 465--471, 1024--1027 of {\it A New Kind of Science.}}
I hope he continues working on this project!

In fact, Wolfram thinks that he not only might be able to find the ultimate TOE,
he might even be able to show that it is the simplest possible TOE!
How does he escape the impact of my results?  Why doesn't Barrow's observation
apply here?
    
First of all, Wolfram is not very interested in proofs, he prefers computational evidence.
Second, Wolfram does not use program-size complexity as his complexity measure.
He uses much more down-to-earth complexity measures.
Third, he is concerned with extremely simple systems, while my methods apply best
to objects with high complexity.
    
Perhaps the best way to explain the difference is to say that he is looking
at ``hardware'' complexity, and I'm looking at ``software'' complexity.  
The objects he studies have complexity less than or equal to that of a universal computer.
Those I study have complexity much larger than a universal computer.
For Wolfram,
a universal computer is the maximum possible complexity, and for me it is the minimum
possible complexity.
    
Anyway, now let's see what's the message from AIT for the working mathematician.

\section*{\bf
VIII Should Mathematics Be More Like Physics? 
Must Mathematical Axioms Be Self-Evident?
}

\footnotesize
``A deep but easily understandable problem about prime
numbers is used in the following to illustrate the parallelism between
the heuristic reasoning of the mathematician and the inductive reasoning
of the physicist\ldots\ [M]athematicians and physicists think alike; they are
led, and sometimes misled, by the same patterns of plausible reasoning.''
    
---George P\'olya, ``Heuristic Reasoning in the Theory
of Numbers'', 1959, reprinted in Alexanderson, {\it The Random Walks of George
P\'olya}, 2000.
    
``The role of heuristic arguments has not been acknowledged
in the philosophy of mathematics, despite the crucial role that they play
in mathematical discovery. The mathematical notion of proof is strikingly
at variance with the notion of proof in other areas\ldots\ Proofs given by
physicists do admit degrees: of two proofs given of the same assertion
of physics, one may be judged to be more correct than the other.''
   
---Gian-Carlo Rota, ``The Phenomenology of Mathematical
Proof'', 1997, reprinted in Jacquette, {\it Philosophy of Mathematics,}
2002, and in Rota, {\it Indiscrete Thoughts,} 1997.
    
``There are two kinds of ways of looking at mathematics\ldots\
the Babylonian tradition and the Greek tradition\ldots\ Euclid discovered that
there was a way in which all the theorems of geometry could be ordered
from a set of axioms that were particularly simple\ldots\ The Babylonian attitude\ldots\
is that you know all of the various theorems and many of the connections
in between, but you have never fully realized that it could all come up
from a bunch of axioms\ldots\ [E]ven in mathematics you can start in different
places\ldots\ In physics we need the Babylonian method, and not the Euclidian
or Greek method.''
    
---Richard Feynman, {\it The Character of Physical Law,}
1965, Chapter 2, ``The Relation of Mathematics to Physics''.
    
``The physicist rightly dreads precise argument, since
an argument which is only convincing if precise loses all its force if
the assumptions upon which it is based are slightly changed, while an argument
which is convincing though imprecise may well be stable under small perturbations
of its underlying axioms.''
    
---Jacob Schwartz, ``The Pernicious Influence of Mathematics
on Science'', 1960, reprinted in Kac, Rota, Schwartz, {\it Discrete Thoughts,}
1992.
    
``It is impossible to discuss realism in logic without
drawing in the empirical sciences\ldots\ A truly realistic mathematics should
be conceived, in line with physics, as a branch of the theoretical construction
of the one real world and should adopt the same sober and cautious attitude
toward hypothetic extensions of its foundation as is exhibited by physics.''
    
---Hermann Weyl, {\it Philosophy of Mathematics and Natural
Science,} 1949, Appendix A, ``Structure of Mathematics'', p.\ 235.
\\

\normalsize
The above quotations are eloquent testimonials to the fact that although
mathematics and physics are different, 
maybe they are not {\bf that} different!
Admittedly,
math organizes our mathematical
experience, which is mental or computational, and physics organizes our physical experience.\footnote
{And in physics everything is an approximation, no equation is exact.}
They are certainly not exactly the same, but maybe it's a matter of degree,
a continuum of possibilities,
and not an absolute, black and white difference.
    
Certainly, as both fields are currently practiced, there is a definite difference in {\bf style}.
But that could change, and is to a certain extent a matter of fashion, not a fundamental
difference.
    
A good source of essays that I---but perhaps not the authors!---regard as generally
supportive of the position that math be considered a branch of physics
is Tymoczko, {\it New Directions
in the Philosophy of Mathematics,} 1998.  In particular
there you will find an essay by Lakatos giving the name ``quasi-empirical''
to this view of the nature of the mathematical enterprise.
    
Why is my position on math ``quasi-empirical''?  Because, as far as I can see,
this is the only way to accommodate the existence
of irreducible mathematical facts gracefully.  Physical postulates are never self-evident,
they are justified pragmatically, and so are close relatives of the not at all
self-evident irreducible
mathematical facts that I exhibited in Section VI. 
    
I'm not proposing that math is a branch of physics just to be controversial.
I was forced to do this against my will!
This happened
in spite of the fact that I'm a mathematician and I love mathematics,
and in spite of the fact that I started with the traditional Platonist position
shared by most working mathematicians.
I'm proposing this because I want mathematics to work better and be more productive.
Proofs are fine, but if you can't find a proof, you should go ahead using heuristic
arguments and conjectures.
    
Wolfram's {\it A New Kind of Science\/} also supports an experimental, quasi-empirical  
way of doing mathematics. This is
partly because Wolfram is a physicist, partly because he believes that unprovable
truths are the rule, not the exception, and 
partly because he believes that our current mathematical theories are highly arbitrary
and contingent.
Indeed, his book may be regarded as a very large
chapter in experimental math. In fact, he had to develop his own programming language,
{\it Mathematica,} to be able to do the massive computations
that led him to his conjectures. 
   
See also Tasi\'c, {\it Mathematics and the Roots of Postmodern Thought,} 2001, for
an interesting perspective on intuition versus formalism.
This is a key question---indeed in my opinion it's an inescapable issue---in any discussion of
how the game of mathematics should be played. 
And it's a question with which I, as a working mathematician, am passionately concerned,
because,
as we discussed in Section VI, formalism has severe limitations.
Only intuition can enable us to go forward and create new ideas and more powerful formalisms.
    
And what are the wellsprings of mathematical intuition and creativity?
In his important forthcoming book on creativity, Tor N{\o}rretranders makes the case that
a peacock, an elegant, graceful woman, and a beautiful mathematical theory,
are all shaped by the same forces, namely what Darwin referred to as ``sexual selection''.
Hopefully this book will be available soon in a language other than Danish!
Meanwhile,
see my dialogue with him in my book {\it Conversations with a Mathematician.}
    
Now, for our last topic, let's look at the entire physical universe!

\section*{\bf
IX Is the Universe Like $\pi$ or Like $\Omega$? Reason versus Randomness!
[Brisson, Meyerstein 1995]
}

\footnotesize
``Parce qu'on manquait d'une d\'efinition rigoreuse
de complexit\'e, celle qu'a propos\'ee la TAI [th\'eorie
algorithmique de l'information], confondre $\pi$
avec $\Omega$ a \'et\'e
plut\^ot la r\`egle que l'exception. Croire, parce que nous avons
ici affaire \`a une croyance, que toutes les suites, puisqu'elles
ne sont que l'encha\^{\i}nement selon une r\`egle rigoureuse de symboles
d\'etermin\'es, peuvent toujours \^etre comprim\'ees
en quelque chose de plus simple, voil\`a la source de l'erreur du
r\'eductionnisme. Admettre la complexit\'e a toujours paru insupportable
aux philosophes, car c'\'etait renoncer \`a trouver un sens rationnel
\`a la vie des hommes.''
    
---Brisson, Meyerstein, {\it Puissance et Limites de la Raison,}
1995, ``Postface. L'erreur du r\'eductionnisme'', p.\ 229.
\\

\normalsize
First let me explain what the number $\Omega$ is.  It's the jewel in AIT's crown,
and it's a number that has attracted a great deal of attention, because
it's a very {\bf dangerous} number!   $\Omega$ is defined to 
be the halting probability of what computer scientists call a universal computer,
or universal Turing machine.\footnote  
{In fact, the precise value of $\Omega$ actually depends on the choice of computer,
and in {\it The Limits of Mathematics\/} I've done that, I've picked one out.}
So $\Omega$ is a probability and therefore it's a real number, 
a number measured with infinite precision,
that's between zero and one.\footnote
{It's ironic that the star of a discrete theory is a real number!
This illustrates the creative tension between the continuous and the discrete.}
That may not sound too dangerous!
    
What's dangerous about $\Omega$ is that (a) it has a simple, straightforward
mathematical definition, but at the same time (b) its numerical value is maximally unknowable,
because a formal mathematical theory whose program-size complexity or explanatory power
is $N$ bits cannot enable you to determine 
more than $N$ bits of the base-two expansion of $\Omega$!
In other words, if you want to calculate $\Omega$, theories
don't help very much, since it takes $N$ bits of theory to get $N$ bits of $\Omega$.
In fact, the base-two bits of $\Omega$ are maximally complex, there's no redundancy,
and $\Omega$ is the
prime example of how unadulterated infinite complexity arises in pure mathematics!
   
How about $\pi$ = 3.1415926\ldots\ the ratio of the circumference of a circle to its diameter?
Well, $\pi$ {\bf looks} pretty complicated, pretty lawless.  For example, all its digits
seem to be equally likely,\footnote 
{In any base all the digits of $\Omega$ are equally likely. 
This is called ``Borel normality''.
For a proof, see my book {\it Exploring Randomness.}
For the latest on $\Omega$, see Calude, {\it Information and Randomness.}}
although this has never been proven.\footnote 
{Amazingly enough, 
there's been some recent progress in this direction by Bailey and Crandall.}
If you are given a bunch of digits from deep inside the decimal
expansion of $\pi$, and you aren't told where they come from,
there doesn't seem to be any
redundancy, any pattern.  But of course, according to AIT, $\pi$ in fact only has
{\bf finite} complexity, because there are algorithms for calculating it with arbitrary
precision.\footnote
{In fact, some terrific new ways to calculate $\pi$  
have been
discovered by Bailey, Borwein and Plouffe. 
$\pi$ lives, it's not a dead subject!}
    
Following Brisson, Meyerstein, {\it Puissance et Limites de la Raison,}
1995, let's now finally discuss whether the physical universe is like $\pi$
= 3.1415926\ldots\ which only has a finite complexity, namely 
the size of the smallest program to generate $\pi$,
or like $\Omega$, which has unadulterated infinite complexity.
Which is it?!
    
Well, if you believe in quantum physics, then Nature plays dice, and that generates
complexity, an infinite amount of it, for example, as frozen accidents, mutations
that are preserved in our DNA.
So at this time most scientists would bet that the universe has infinite complexity,
like $\Omega$ does.
But then the world is incomprehensible, or at least a large part of it will always
remain so, the accidental part, all those frozen accidents, the contingent part.
    
But some people still hope that the world has finite complexity like $\pi$, it just
{\bf looks} like it has high complexity.  
If so, then we might eventually be able to comprehend everything, and there is an ultimate TOE!
But then you have to believe that quantum mechanics is wrong, as currently practiced,
and that all that quantum randomness is really only {\bf pseudo-randomness}, like 
what you find in
the digits of $\pi$. You have to believe that the world is actually deterministic,
even though our current scientific theories say that it isn't!
    
I think Vienna physicist Karl Svozil feels that way [private communication; see his
{\it Randomness \& Undecidability in Physics,} 1994].  I know
Stephen Wolfram does, he says so in his book.  Just take a look at the discussion of
fluid turbulence and of the second law of thermodynamics in {\it A New Kind of Science.}
Wolfram believes that very simple deterministic algorithms ultimately account for all
the apparent complexity we see around us, just like they do in $\pi$.\footnote 
{In fact, Wolfram himself explicitly makes the connection with $\pi$.
See {\bf meaning of the universe} on page 1027 of {\it A New Kind of Science.}}
He believes that the world {\bf looks}
very complicated, but is actually very simple. 
There's no randomness, there's only pseudo-randomness.
Then nothing is contingent, everything
is necessary, everything happens for a reason. [Leibniz!] 
   
Who knows!  Time will tell! 
   
Or perhaps from {\bf inside} this world we will never be able to tell the difference,
only an {\bf outside} observer could do that [Svozil, private communication].

\section*{\bf
Postscript
}

Readers of this paper may enjoy the somewhat different perspective in 
my chapter ``Complexit\'e, logique et hasard'' 
in Benkirane, {\it La Complexit\'e.} 
Leibniz is there too.
    
In addition, see my {\it Conversations with a Mathematician,} a book on
philosophy disguised as a series of dialogues---not the first time that this has happened!

Last but not least, see Zwirn, {\it Les Limites de la Connaissance,} 
that also supports the thesis that understanding is compression,
and the masterful multi-author two-volume
work, {\it Kurt G\"odel, Wahrheit \& Beweisbarkeit,} 
a treasure trove of information about G\"odel's life and work.

\section*{\bf
Acknowledgement
}

Thanks to Tor N{\o}rretranders for providing the original German for the Einstein
quotation at the beginning of this paper, and also the word for word translation.

The author is grateful to Fran\c{c}oise Chaitin-Chatelin for innumerable stimulating
philosophical discussions. He dedicates this paper to her unending quest
to understand.

\section*{\bf
Bibliography
}

\footnotesize
\begin{itemize}
\item
Gerald W. Alexanderson, {\it The Random Walks of George P\'olya,}
MAA, 2000.
\item  
John D. Barrow, Frank J. Tipler, {\it The Anthropic Cosmological Principle,} Oxford University
Press, 1986.
\item  
John D. Barrow, {\it Between Inner Space and Outer Space,} Oxford University
Press, 1999.
\item  
Reda Benkirane, {\it La Complexit\'e, Vertiges et Promesses,} Le Pommier, 2002.
\item  
Max Born, {\it Experiment and Theory in Physics,} Cambridge University Press, 1943.
Reprinted by Dover, 1956.
\item  
Luc Brisson, F. Walter Meyerstein, {\it Inventer l'Univers,} Les Belles
Lettres, 1991.
\item  
Luc Brisson, F. Walter Meyerstein, {\it Inventing the Universe,} SUNY
Press, 1995.
\item  
Luc Brisson, F. Walter Meyerstein, {\it Puissance et Limites de la Raison,}
Les Belles Lettres, 1995.
\item  
F. Br\'ody, T. V\'amos, {\it The Neumann Compendium,} World Scientific, 1995.
\item  
Bernd Buldt et al., {\it Kurt G\"odel, Wahrheit \& Beweisbarkeit. Band 2: Kompendium
zum Werk,} \"obv \& hpt, 2002.
\item  
Cristian S. Calude, {\it Information and Randomness,} Springer-Verlag, 2002.
\item  
Gregory J. Chaitin, {\it The Limits of Mathematics, The Unknowable,
Exploring Randomness, Conversations with a Mathematician,} Springer-Verlag,
1998, 1999, 2001, 2002.
\item  
John Cornwell, {\it Nature's Imagination,} Oxford University Press,
1995.
\item  
COSRIMS, {\it The Mathematical Sciences,} MIT Press, 1969.
\item  
Albert Einstein, {\it Ideas and Opinions,} Crown, 1954. Reprinted by Modern Library, 1994.
\item  
Albert Einstein, {\it Autobiographical Notes,} Open Court, 1979.
\item  
Richard Feynman, {\it The Character of Physical Law,} MIT Press, 1965.
Reprinted by Modern Library, 1994, with a thoughtful introduction by James Gleick. 
\item  
Richard P. Feynman, Robert B. Leighton, Matthew Sands, {\it The Feynman
Lectures on Physics,} Addison-Wesley, 1963.
\item  
Dale Jacquette, {\it Philosophy of Mathematics,} Blackwell, 2002.
\item  
Mark Kac, Gian-Carlo Rota, Jacob T. Schwartz, {\it Discrete Thoughts,} Birkh\"auser, 1992.
\item  
Eckehart K\"ohler et al., {\it Kurt G\"odel, Wahrheit \& Beweisbarkeit. Band 1: Dokumente und
historische Analysen,} \"obv \& hpt, 2002.
\item  
Bernd-Olaf K\"uppers, {\it Information and the Origin of Life,} MIT Press, 1990. 
\item  
G. W. Leibniz, {\it Philosophical Essays,} edited and translated by
Roger Ariew and Daniel Garber, Hackett, 1989.
\item  
Ernst Mach, {\it The Science of Mechanics,} Open Court, 1893.
\item
Paolo Mancosu, {\it From Brouwer to Hilbert,} Oxford University Press, 1998.
\item
Karl Menger, {\it Reminiscences of the Vienna Circle and the Mathematical Colloquium,} Kluwer, 1994.
\item  
James R. Newman, {\it The World of Mathematics,} Simon and Schuster, 1956.
Reprinted by Dover, 2000.
\item  
Karl R. Popper, {\it The Logic of Scientific Discovery,} Hutchinson Education, 1959. 
Reprinted by Routledge, 1992.
\item  
Gian-Carlo Rota, {\it Indiscrete Thoughts,} Birkh\"auser, 1997.
\item  
Paul Arthur Schilpp, {\it Albert Einstein, Philosopher-Scientist,} Open Court, 1949.
\item  
Karl Svozil, {\it Randomness \& Undecidability in Physics,} World Scientific, 1994.
\item  
Vladimir Tasi\'c, {\it Mathematics and the Roots of Postmodern Thought,}
Oxford University Press, 2001.
\item  
Thomas Tymoczko, {\it New Directions in the Philosophy of Mathematics,}
Princeton University Press, 1998.
\item  
Hermann Weyl, {\it The Open World,} Yale University Press, 1932. 
Reprinted by Ox Bow Press, 1989.
\item  
Hermann Weyl, {\it Philosophy of Mathematics and Natural Science,} 
Princeton University Press, 1949.
\item  
Stephen Wolfram, {\it A New Kind of Science,} Wolfram Media, 2002.
\item  
Herv\'e Zwirn, {\it Les Limites de la Connaissance,} Odile Jacob, 2000. 
\end{itemize}

\end{document}